\documentclass[10pt]{amsart}

\usepackage{latexsym, amsmath,amssymb}

\newcommand{\cal}{\mathcal}

\renewcommand{\epsilon}{\varepsilon}
\newcommand{\newsection}[1]
{\subsection{#1}\setcounter{theorem}{0} \setcounter{equation}{0}
\par\noindent}

\newtheorem{theorem}{Theorem}

\newtheorem{lemma}[theorem]{Lemma}
\newtheorem{corr}[theorem]{Corollary}

\newtheorem{proposition}[theorem]{Proposition}
\newtheorem{deff}[theorem]{Definition}

\newcommand{\bth}{\begin{theorem}}
\newcommand{\ble}{\begin{lemma}}
\newcommand{\bcor}{\begin{corr}}

\newcommand{\bdeff}{\begin{deff}}

\newcommand{\bprop}{\begin{proposition}}
\newcommand{\ele}{\end{lemma}}
\newcommand{\ecor}{\end{corr}}
\newcommand{\edeff}{\end{deff}}

\newcommand{\eprop}{\end{proposition}}

\newcommand{\cd}{\, \cdot\, }

\newcommand{\Z}{{\mathbb Z}}
\newcommand{\Rn}{{\mathbb R}^n}

\newcommand{\supp}{\text{supp}}
\renewcommand{\Pi}{\varPi}

\renewcommand{\epsilon}{\varepsilon}

\newcommand{\R}{{\mathbb R}}
\newcommand{\g}{{\mathrm g}}

\begin{document}

\title[Abstract Strichartz Estimates]
{On Abstract Strichartz Estimates and the Strauss Conjecture for 
Nontrapping Obstacles}

\thanks{The first author was supported in part by
the Grant-in-Aid for Young Scientists (B) (No.18740069),
The Ministry of Education, Culture,
Sports, Science and Technology, Japan, and he would like to thank
the Department of Mathematics at the Johns Hopkins University
for the hospitality and financial support during his visit
where part of this research was carried out.  The second, third and fourth 
authors were supported by the National Science Foundation.  The fifth author
was supported by project 10728101 of NSFC and the ``111''
project and Doctoral Programme Foundation of the Ministry of Education of 
China.}
\author[K. Hidano]{Kunio Hidano}
\address{Department of Mathematics, Mie University}
\author[J. Metcalfe]{Jason Metcalfe}
\address{Department of Mathematics, University of North Carolina, Chapel Hill}
\author[H. F. Smith]{Hart F. Smith}
\address{Department of Mathematics, University of Washington, Seattle}
\author[C. D. Sogge]{Christopher D. Sogge}
\address{Department of Mathematics, Johns Hopkins University}
\author[Y. Zhou]{Yi Zhou}
\address{School of Mathematical Science, Fudan University}
%\authorrunning{K. Hidano, J. Metcalfe, H. F. Smith, C. D. Sogge and Y. Zhou}

\maketitle
\newsection{Introduction}

The purpose of this paper is to show how local energy decay estimates for 
certain linear wave equations involving compact perturbations of the standard 
Laplacian lead to optimal global existence theorems for the corresponding 
small amplitude nonlinear wave equations with power nonlinearities.  To 
achieve this goal, at least for spatial dimensions $n=3$ and $4$, we shall 
show how the aforementioned linear decay estimates can be combined with 
``abstract Strichart'' estimates for the free wave equation to prove 
corresponding estimates for the perturbed wave equation when $n\ge3$.  
As we shall see, we are only partially successful in the latter endeavor 
when the dimension is equal to two, and therefore, at present, our 
applications to nonlinear wave equations in this case are limited.

Let us start by describing the local energy decay assumption that we 
shall make throughout.  We shall consider wave equations on the exterior 
domain $\Omega\subset \R^n$ of a compact obstacle:
\begin{equation}\label{1}
\begin{cases}
(\partial^2_t-\Delta_\g)u(t,x)=F(t,x), \quad (t,x)\in \R_+\times \Omega
\\
u(0,\cd)=f %\in \dot H^\gamma_D(\Omega)
\\
\partial_t u(0,\cd)=g%\in \dot H^{\gamma-1}_D(\Omega)
\\
(Bu)(t,x)=0,\quad \text{on } \, \R_+\times \partial \Omega,
\end{cases}
\end{equation}
where for simplicity we take $B$ to either be the identity operator 
(Dirichlet-wave equation) or the inward pointing normal derivative
$\partial_\nu$ (Neumann-wave equation).  We shall also assume throughout 
that the spatial dimension satisfies $n\ge 2$.

The operator $\Delta_\g$ is the Laplace-Beltrami operator associated with a 
smooth, time independent Riemannian metric ${\mathbf g}_{jk}(x)$ which we 
assume equals the Euclidean metric $\delta_{jk}$ for $|x|\ge R$, some $R$.  
The set $\Omega$ is assumed to be either all of $\Rn$, or else 
$\Omega=\Rn\backslash {\mathcal K}$ where ${\mathcal K}$ is a compact 
subset of $|x|<R$ with smooth boundary.

We can now state the main assumption that we shall make.

\noindent{\bf Hypothesis B.} {\it Fix the boundary operator} $B$ 
{\it and the exterior domain} $\Omega \subset \R^n$ {\it as above.  
We then assume that given} $R_0>0$
\begin{multline}\label{1.2}
\int_0^\infty \Bigl(\, \|u(t,\cd)\|^2_{H^1(|x|<R_0)}+
\|\partial_tu(t,\cd)\|^2_{L^2(|x|<R_0)}\Bigr)\, dt
\\
 \lesssim \|f\|^2_{H^1} + \|g\|^2_{L^2} +\int_0^\infty\|F(s,\cd)\|_{L^2}^2\, 
ds,
\end{multline}
{\it whenever} $u$ {\it is a solution of} \eqref{1} {\it with data} 
$(f(x),g(x))$ {\it and forcing term} $F(t,x)$ {\it that both vanish for} 
$|x|>R_0$.

Here $A\lesssim B$ means that $A$ is bounded by a constant times $B$, and, 
in what follows, the constant might change at each occurrence.  Also, 
$\|h\|_{H^1(|x|<R_0)}$ denotes the $L^2$-norm of $h$ and $\nabla_x h$ over the 
set $\{x\in \Omega: \, |x|<R_0\}$.

Let us review some important cases where the assumption \eqref{1.2} is valid.  
First of all, results from Vainberg \cite{V}, combined with the propagation 
of singularity results of Melrose and Sj\"ostrand \cite{MelSj}, imply that if 
$\Delta$ is the standard Euclidean Laplacian and $\Omega$ is nontrapping, then 
if $u$ is a solution of \eqref{1} with data of fixed compact support and 
forcing term $F\equiv 0$, then with $u'=(\partial_tu,\nabla_xu)$,
$$\|u'(t,\cd)\|_{L^2(|x|<R_0)}\le \alpha(t) \|u'(0,\cd)\|_{L^2},$$
where $\alpha(t)=O((1+t)^{-(n-1)})$ for either the Dirichlet-wave equation 
or the Neumann-wave equation when $n\ge3$.  For $n=2$, if 
$\partial\Omega$ is assumed to be nonempty one has 
$\alpha(t)=O((\log(2+t))^{-2}(1+t)^{-1})$ for the Dirichlet-wave equation.
Here we have used that, due to the Dirichlet boundary conditions and the 
fundamental theorem
of calculus, the local $L^2$ norm can be controlled by the local $L^2$ 
norm of the gradient.  
Since these bounds yield $\alpha(t)\in L^1(\R_+)$, we conclude that 
Hypothesis B is valid in these cases.  We remark that when $\Omega=\R^2$  
and $\Delta_\g=\Delta$, then
$\alpha \approx t^{-1}$ for large $t$ (see \cite{Ral}), 
and so, in this case, $\alpha\notin L^1(\R_+)$.  Proofs of these results for 
$n\ge 3$ can be found in Melrose \cite{Mel} and Ralston \cite{Ral}, while the 
result for the Dirichlet-wave equation for $n=2$ follows from 
Vainberg \cite{V} (see \S 4 and Remark 4 on p. 40). \footnote{
We are very grateful to Jim Ralston for patiently explaining these results 
and their history to us.}

For the case where $\Delta_\g$ is assumed to be a time-independent 
variable coefficient compact perturbation of $\Delta$ and $\Omega$ is 
assumed to be nontrapping with respect to the metric associated with 
$\Delta_\g$, one also has that \eqref{1.2} is valid for the Dirichlet-wave 
equation for all $n\ge3$ as well for $n=2$ if
$\partial \Omega\ne \emptyset$.  See Taylor \cite{taylor} and
Burq \cite{burq}.

Having described the main assumption about the linear problem, let us now 
describe the nonlinear equations that we shall consider.  They are of the form
\begin{equation}\label{1.3}
\begin{cases}
(\partial_t^2-\Delta_\g)u(t,x)=F_p\bigl(u(t,x)\bigr), 
\quad (t,x)\in \R_+\times \Omega
\\
Bu=0, \quad \text{on } \, \R_+\times \partial\Omega
\\
u(0,x)=f(x), \quad \partial_t u(0,x)=g(x), \quad x\in \Omega,
\end{cases}
\end{equation}
with $B$ as above.  We shall assume that the nonlinear term behaves like 
$|u|^p$ when $u$ is small, and so we assume that
\begin{equation}\label{1.4}
\sum_{0\le j\le 2} |u|^j\, \bigl|\, \partial^j_u F_p(u)\, \bigr| \,  
\lesssim \, |u|^p
\end{equation}
when $u$ is small.  %We remark that if $n=3$, the proofs to follow show that 
%this assumption can be weakened in the existence results to be $F_p(0)=0$ 
%and $F_p'(u)=O(|u|^{p-1})$.

We shall be assuming that the data (and some of its derivatives) are small 
in certain Sobolev norms that we now describe.  

As in the earlier works that proved global Strichartz estimates 
(\cite{burq},  \cite{M}, \cite{SS}), we shall restrict ourselves to the 
case where the Sobolev index $\gamma$ is smaller than $n/2$.  One reason for 
this is that the Strichartz estimates that seem to arise in applications 
always have $\gamma\le 1$.  Another reason is that when $|\gamma|<n/2$ , 
multiplication 
by a smooth function $\beta \in C^\infty_0(\Rn)$ is continuous from 
$\dot H^\gamma(\Rn)$ to $H^\gamma(\Rn)$ and the two norms are equivalent 
on functions with fixed compact support.  Recall that $\dot H^\gamma(\Rn)$ 
is the homogeneous Sobolev space with norm given by
$$\|f\|^2_{\dot H^\gamma(\Rn)}=\bigl\|\, 
(\sqrt{-\Delta})^{\gamma} f\bigr\|_{L^2(\Rn)}^2
=(2\pi)^{-n}\int_{\Rn}\bigl| \, |\xi|^{\gamma}\Hat f(\xi)\, \bigr|^2\, d\xi,$$
while the inhomogeneous Sobolev space $H^\gamma(\Rn)$ has norm defined by
$$\|f\|_{H^\gamma(\Rn)}^2 = \bigl\| \, (1-\Delta)^{\gamma/2}f\, 
\bigr\|^2_{L^2(\Rn)}
= (2\pi)^{-n}\int_{\Rn} \bigl|\, (1+|\xi|^2)^{\gamma/2}\hat f(\xi)\, 
\bigr|^2\, d\xi,$$
with $\Hat f$ denoting the Fourier transform and $\Delta$ denoting 
the standard Laplacian.

Let us now describe the Sobolev spaces on $\Omega$ that
we shall consider. Let $\beta$ be
a smooth cutoff on $\R^n$ with $\beta$ and $1-\beta$ respectively supported
where $|x|<2R$ and $|x|>R$. Let $\Omega'$ be the embedding
of $\Omega\cap\{|x|<2R\}$ into the torus obtained by
periodic extension of $\Omega\cap [-2R,2R]^n$, so that 
$\partial\Omega' = \partial\Omega$. We define
\begin{align*}
\|f\|_{H_B^\gamma(\Omega)}&=\|\beta f\|_{H_B^\gamma(\Omega')}
+\|(1-\beta)f\|_{H^\gamma(\R^n)}\\
\|f\|_{\dot{H}_B^\gamma(\Omega)}&=\|\beta f\|_{H_B^\gamma(\Omega')}
+\|(1-\beta)f\|_{\dot{H}^\gamma(\R^n)}\,,\qquad |\gamma|< n/2\,.
\end{align*}
The spaces $H^\gamma_B(\Omega')$ are defined by a spectral decomposition
of $\Delta_\g|_{\Omega'}$ subject to the boundary condition $B$.
In the homogeneous spaces $\dot{H}_B^\gamma(\Omega)$ it is assumed that
$(1-\beta)f$ belongs to $\dot{H}^\gamma(\R^n)$, so that the Sobolev
embedding $\dot{H}^\gamma_B(\Omega)\hookrightarrow L^p(\Omega)$ holds with
$p=2n/(n-2\gamma)$.
From this, it is verified that the Sobolev spaces on $\Omega$
are independent of the
choice of $\beta$ and $R$, and thus the $\dot{H}^\gamma_B(\Omega)$ and
$H^\gamma_B(\Omega)$
norms are equivalent on functions of fixed bounded support.
We note that $H^{-\gamma}_B(\Omega)$ is the dual of $H^\gamma_B(\Omega)$, and
$\dot{H}^{-\gamma}_B(\Omega)$ is dual to $\dot{H}^\gamma_B(\Omega)$ for
$|\gamma|<n/2$. Also, for $\gamma$ a nonnegative integer,
\begin{align*}
\|f\|_{H_B^\gamma(\Omega)}^2&\approx 
\sum_{|\alpha|\le\gamma}\|\partial^\alpha_x\! f\|_{L^2(\Omega)}^2
\\
\|f\|_{\dot{H}_B^\gamma(\Omega)}^2&\approx
\sum_{|\alpha|=\gamma}\|\partial^\alpha_x\! f\|_{L^2(\Omega)}^2.
\end{align*}
%[H] addition follows
The Sobolev spaces as defined are verified to be an analytic interpolation
scale of spaces. The above definition then agrees, for nonnegative
integer $\gamma$, with the subspace of $H^\gamma(\overline\Omega)$
such that $B(\Delta_\g^j f)=0$ for all $j$ for which the trace
is well defined, and for general $\gamma$ by duality and interpolation.
Finally,  for every $\gamma$ the set of functions 
$f\in C_0^\infty(\overline{\Omega})$
such that $B(\Delta_\g^j f)=0$ for all $j\ge 0$ is dense in the norm.

Our hypotheses regarding the data in \eqref{1.3} will only involve 
certain $\gamma\in(0,\frac 12)$, while the ones in the 
abstract Strichartz estimates 
to follow only involve certain $\gamma \le(n-1)/2$.  In practice the useful 
Strichartz-type estimates always involve $\gamma\in (0,1]$.  This is the case 
for the mixed-norm Strichartz estimates of Keel and Tao \cite{KT} and others 
for the case
$\Omega=\Rn$, $\Delta=\Delta_\g$, as well as for the mixed-norm estimates 
for \eqref{1} that we shall state.

The data $(f,g)$ in Theorem \ref{exist} below
will have second derivatives belonging
to $\dot{H}^\gamma_B(\Omega)\times\dot{H}^{\gamma-1}_B(\Omega)$, 
where $\gamma\in(0,\frac 12)$, thus will locally belong to 
$H^{2+\gamma}(\overline{\Omega})\times H^{1+\gamma}(\overline{\Omega})$. 
The boundary condition for $(f,g)$ to locally belong to 
$H^{2+\gamma}_B(\Omega)\times H^{1+\gamma}_B(\Omega)$ for 
$\gamma\in (0,\frac 12)$ 
is the same as for 
$H^2_B(\Omega)\times H^1_B(\Omega)$, which for the Dirichlet case
is $f|_{\partial\Omega}=g|_{\partial\Omega}=0$, and for Neumann is
$\partial_\nu f|_{\partial\Omega}=0$. These are the assumptions placed
on the data $(f,g)$ in Theorem \ref{exist}.

If we let 
$$\{Z\} = \{ \partial_l, \, x_j\partial_k-x_k\partial_j: 
\, 1\le l\le n, \, 1\le j<k\le n\}$$
then we can now state our existence theorem for \eqref{1.3}.

\begin{theorem}\label{exist}  Let $n=3$ or $4$, and fix $\Omega\subset \Rn$ 
and boundary operator $B$ as above.  Assume further that Hypothesis B is valid.

Let $p=p_c$ be the positive root of 
\begin{equation}\label{1.5}
(n-1)p^2-(n+1)p-2=0,
\end{equation}
and fix $p_c<p< (n+3)/(n-1)$.
Then if
\begin{equation}\label{1.6}
\gamma=\tfrac{n}2-\tfrac2{p-1},
\end{equation}
there is an $\varepsilon_0>0$ depending on $\Omega, B$ and $p$ so that 
\eqref{1.3} has a global solution satisfying 
$(Z^\alpha u(t,\cd), \partial_t Z^\alpha u(t,\cd))\in 
\dot H^\gamma_B\times \dot H^{\gamma-1}_B$, $|\alpha|\le 2$, $t\in \R_+$, 
whenever the initial data satisfies the boundary conditions of order $2$,
and
\begin{equation}\label{1.7}
\sum_{|\alpha|\le 2}\Bigl(\, \|Z^\alpha f\|_{\dot H^\gamma_B(\Omega)}
+\|Z^\alpha g\|_{\dot H^{\gamma-1}_B(\Omega)}\, \Bigr)<\varepsilon
\end{equation}
with $0<\varepsilon<\varepsilon_0$.  \end{theorem}

In the case where $\Omega=\Rn$ and $\Delta_\g=\Delta$ it is known that $p>p_c$ 
is necessary for global existence (see John \cite{John}, Glassey \cite{Glassey}
and Sideris \cite{Sideris}).  In this case under a somewhat more restrictive 
smallness condition global existence was established by John \cite{John} for 
the case where $n=3$, then Glassey \cite{Glassey} for $n=2$, Zhou \cite{Zhou} 
for $n=4$, Lindblad and Sogge \cite{LS2} for $n\le 8$ and then Georgiev, 
Lindblad and Sogge \cite{GLS} for all $n$ (see also Tataru \cite{Ta}).  
For obstacles, when $n=4$, $\Delta_\g=\Delta$ the results in Theorem 
\ref{exist} for the Dirichlet-wave equation outside of nontrapping obstacles 
under a somewhat more restrictive smallness assumption was obtained in 
\cite{MSYY}.

Also, when $\Omega=\R^3$, $\Delta_\g=\Delta$, it was shown in Sogge \cite{So} 
that, for the spherically symmetric case, the variant of the condition 
\eqref{1.7} saying that the 
$\dot H^\gamma(\R^3)\times \dot H^{\gamma-1}(\R^3)$ norm of the data be 
small with $\gamma$ as in \eqref{1.6} is sharp.  Further work in this 
direction (for the non-obstacle case) was done by Hidano \cite{H}, \cite{H2} 
and Fang and Wang \cite{FW}.

It is not difficult 
to see that the condition \eqref{1.7} is sharp in the sense that there are 
no global existence results for $\gamma>\tfrac{n}2-\tfrac2{p-1}$.  
To do this we use well known results concerning blowup solutions for 
$(\partial_t^2-\Delta)v=|v|^p$, $p>0$, in $\R_+\times \R^n$ (see Levine 
\cite{levine}).  Specifically, we shall use the fact that given $\delta>0$ 
one can find $C^\infty_0$ data $(v_0,v_1)$ vanishing for $|x|<R$ so that the 
solution of $(\partial^2_t-\Delta)v=|v|^p$, 
$v(0,\cd)=v_0, \partial_tv(0,\cd)=v_1$ blows up within time $\delta$.
Next, let us assume that the above global existence results for 
$(\Omega,B,\Delta_\g)$ held for this nonlinearity and some 
$\gamma>\tfrac{n}2-\tfrac{2}{p-1}$ in \eqref{1.7}.  Then, if $\lambda$ is 
sufficiently large, the $\dot H^\gamma_B \times \dot H^{\gamma-1}_B$ norm of 
$(\lambda^{-2/(p-1)}v_0(\cd/\lambda), \, \lambda^{-1-2/(p-1)}v_1(\cd/\lambda))$
 would be bounded by its $\dot H^\gamma(\Rn)\times \dot H^{\gamma-1}(\Rn)$ norm,
which equals $\lambda^{n/2-2/(p-1)-\gamma}\|(v_0,v_1)\|_{\dot H^\gamma(\R^n)
\times \dot H^{\gamma-1}(\Rn)} $.
 Since this goes to zero as $\lambda\to \infty$ for 
$\gamma >\tfrac{n}2-\tfrac{2}{p-1}$, we conclude that if the above existence 
results held for this value of $\gamma$ then we would obtain a global solution 
of $(\partial_t^2-\Delta_\g)u_\lambda=|u_\lambda|^p$, 
$u_\lambda(t,x)=0, (t,x)\in \R_+\times\partial\Omega$  with initial data 
$(\lambda^{-2/(p-1)}v_0(\cd/\lambda), \lambda^{-1-2/(p-1)}v_1(\cd/\lambda))$.  
Since $v_0$ and $v_1$ vanish for $|x|<R$, by finite propagation speed, if 
$\delta>0$ is small and fixed, then for large $\lambda$ if we extend 
$u_\lambda$ to be zero on $\Omega^c$ then the resulting function would agree 
with the solution of the Minkowski space wave equation 
$(\partial_t^2-\Delta)v_\lambda = |v_\lambda|^p$ on 
$[0,\delta\lambda]\times \R^n$ with data 
$(\lambda^{-2/(p-1)}v_0(\cd/\lambda), \lambda^{-1-2/(p-1)}v_1(\cd/\lambda))$.  
By scaling $v(t,x)=\lambda^{2/(p-1)}v_\lambda(\lambda t, \lambda x)$ would 
then solve the Minkowski space equation $(\partial_t-\Delta)v=|v|^p$ on 
$[0,\delta]\times \Rn$ with initial data $(v_0(x),v_1(x))$.  
As we noted before, 
we can always choose $(v_0,v_1)$ so that this is impossible for a given 
$\delta>0$, which allows us to conclude that the above existence results 
do not hold if the Sobolev exponent $\gamma$ in \eqref{1.7} is larger than 
$\tfrac{n}2-\tfrac2{p-1}$.

As a final remark, we point out that we have restricted ourselves to the 
case where $p<(n+3)/(n-1)$ because of the techniques that we shall employ.  
However, since the solutions obtained are small, the above existence theorem 
leads to small-data global existence of \eqref{1.3} when $p$ is larger than 
or equal to the conformal power $(n+3)/(n-1)$.

\medskip

To prove Theorem \ref{exist}, we shall use certain ``abstract Strichartz 
estimates'' which we now describe.  Earlier works (\cite{burq}, \cite{M}, 
\cite{SS})  have focused on establishing certain mixed norm, $L^q_tL^r_x$ 
estimates on $\R_+\times \Omega$ for solutions of \eqref{1}.  For certain 
applications, such as obtaining the Strauss conjecture in various settings, 
it is convenient to replace the $L^r_x$ norm with a more general one.  
To this end, we consider pairs of normed function spaces $X(\R^n)$
and $X(\Omega)$. The spaces are localizable,
in that $\|f\|_X\approx \|\beta f\|_X+\|(1-\beta)f\|_X$ for smooth, compactly
supported $\beta$, with $\beta=1$ on a neighborhood of 
$\R^n\backslash\Omega$ in case $X=X(\Omega)$.
Finally, we assume that
\begin{equation}\label{local}
\|(1-\beta)f\|_{X(\Omega)}\approx\|(1-\beta)f\|_{X(\R^n)}
\end{equation}
for such $\beta$. Weighted mixed $L^p$ spaces, as well as
$\bigl(\dot H^\gamma(\R^n),\dot H^\gamma_B(\Omega)\bigr)$, are the examples
used in the proof of Theorem \ref{exist}.

We shall let $\|\cdot\|_{X'}$ denote the dual norm (respectively over
$\R^n$ and $\Omega$) so that
$$\|u\|_X = \sup_{\|v\|_{X'}=1} \Bigl| \, 
\int u\, \overline{v}\, dx\, \Bigr|\,.$$
An important example for us is when 
$$\|u\|_{X}=\|\, |x|^\alpha u\|_{L^p},$$
for a given $1\le p\le \infty$ and $|\alpha|<n/p$, in which case the dual 
norm is
$$\|v\|_{X'}=\| \, |x|^{-\alpha} v\|_{L^{p'}},$$
with $p'$ denoting the conjugate exponent.

We shall consider time Lebesgue exponents $q\ge2$ and assume that we have 
the global Minkowski abstract Strichartz estimates
\begin{equation}\label{1.8}
\|v\|_{L^q_tX(\R\times \Rn)}\lesssim \|v(0,\cd)\|_{\dot H^{\gamma}(\Rn)}
+\|\partial_tv(0,\cd)\|_{\dot H^{\gamma-1}(\Rn)}\,,
\end{equation}
assuming that
\begin{equation}\label{1.9}(\partial_t^2-\Delta)v=0 \quad \text{in} \;\; 
\R\times \R^n\,.
\end{equation}
Here
$$
\|v\|_{L^q_tX(I\times \Rn)}= \Bigl(\, \int_I \|v(t,\cd)\|^q_X\, dt \, 
\Bigr)^{1/q}, \quad I\subset \R.
$$
We shall also consider analogous norms on $I\times \Omega$, $I\subset \R$,
$$\|u\|_{L^q_tX(I\times \Omega)}=
\Bigl(\, \int_I\|u(t,\cd)\|^q_{X(\Omega)}\, dt\, \Bigr)^{1/q}.$$

In addition to Hypothesis B 
and \eqref{1.8}, we shall assume that we have the local abstract 
Strichartz estimates for $\Omega$:
\begin{equation}\label{1.10}
\|u\|_{L^q_tX([0,1]\times \Omega)}\lesssim \|f\|_{\dot H^\gamma_B(\Omega)}
+\|g\|_{\dot H^{\gamma-1}_B(\Omega)},
\end{equation}
assuming that $u$ solves \eqref{1} with vanishing forcing term, i.e.,
\begin{equation}\label{1.11}
(\partial^2_t-\Delta_\g)u=0\quad \text{in } \, \, [0,1]\times \Omega.
\end{equation}
 
\begin{deff}\label{def1.2} 
When \eqref{1.8} and \eqref{1.10} hold we say that 
$(X,\gamma,q)$ is an admissible triple.
\end{deff}

We can now state our main estimate.

\begin{theorem}\label{theorem1.3}  Let $n\ge2$ and assume that $(X,\gamma,q)$ 
is an admissible triple with 
\begin{equation}\label{1.12} q>2\quad \text{and} \quad 
\gamma\in[-\tfrac{n-3}2,\tfrac{n-1}2].
\end{equation}
Then if Hypothesis B is valid and if $u$ solves \eqref{1} with 
$(\partial^2_t-\Delta_\g)u\equiv 0$,
we have the global abstract Strichartz estimates 
\begin{equation}\label{1.13}
\|u\|_{L^q_tX(\R\times \Omega)}\lesssim \|f\|_{\dot H^\gamma_B(\Omega)}+
\|g\|_{\dot H^{\gamma-1}_B(\Omega)}.
\end{equation}
\end{theorem}

The condition on $\gamma$ in \eqref{1.12} is the one to ensure that $\gamma$ 
and $1-\gamma$ are both $\le (n-1)/2$, which is what the proof seems to 
require.  Unfortunately, for $n=2$, this forces $\gamma$ to be equal to $1/2$, 
while a larger range of $\gamma\in (0,1)$ is what certain applications require.
For this reason, we are unable at present to show that the Strauss conjecture 
for obstacles holds when $n=2$.   See the end of the next section for 
further discussion.

\begin{corr}\label{corollary1.4}  Assume that $(X,\gamma,q)$ and 
$(Y,1-\gamma,r)$ are admissible triples and that Hypothesis B is valid.  
Also assume that \eqref{1.13} holds for $(X,\gamma,q)$ and $(Y,1-\gamma,r)$,
and that $0\le \gamma\le 1$.
Then we have the following 
global abstract Strichartz estimates for the solution of \eqref{1}
\begin{equation}\label{1.15}
\|u\|_{L^q_tX(\R_+\times\Omega)}\lesssim 
\|f\|_{\dot H^{\gamma}_B(\Omega)}
+\|g\|_{\dot H^{\gamma-1}_B(\Omega)}
+\|F\|_{L^{r'}_tY'(\R_+\times \Omega)},
\end{equation}
where $r'$ denotes the conjugate exponent to $r$ and 
$\|\, \cd\|_{Y'}$ is the dual norm to $\|\cd\|_Y$.
\end{corr}

%A special case that of mixed-norms which were treated in  \cite{burq},  
%\cite{M} and \cite{SS}, where one takes $X=L^r(\Rn)$ and $Y=L^s(\Rn)$ for 
%appropriate values of $r$ and $s$.

%Note that by causality, the estimate also holds with $\R$ replaced by $\R_+$.
For simplicity, in the corollary we have limited ourselves to the case where 
$0\le \gamma\le 1$ since that is all that is needed for the applications.  
%By doing this the only compatibility condition on the data is that
%$f(x)=0, \, x\in \partial \Omega$ if $\gamma_X>1/2$ when $B=Id$; no condition 
%is needed for the Neumann case, $B=\partial_\nu$.

Let us give the simple argument that shows that \eqref{1.15} follows from 
\eqref{1.13}. 
To prove \eqref{1.15}, we may assume by \eqref{1.13} 
that the initial data vanishes.
By \eqref{1.13} and the Duhamel formula, if $P=\sqrt{-\Delta_\g}$ is
the square root of minus the Laplacian (with the boundary conditions $B$), 
then we need show
$$
\Bigl\|\,\int_0^t \sin\bigl((t-s)P\bigr)P^{-1}F(s,\cd)\, ds\, 
\Bigr\|_{L^{q}_tX(\R_+\times\Omega)}\lesssim
\|F\|_{L^{r'}_tY'(\R_+\times \Omega)}\,.
$$
Since $q>r'$, an application of the Christ-Kiselev lemma 
(cf. \cite{CK}, \cite{SS}, \cite[chapter 4]{So2}) shows that it
suffices to prove the estimate
$$
\Bigl\|\,\int_0^\infty \sin\bigl((t-s)P\bigr)P^{-1}F(s,\cd)\, ds\, 
\Bigr\|_{L^{q}_tX(\R_+\times\Omega)}\lesssim
\|F\|_{L^{r'}_tY'(\R_+\times \Omega)}\,.
$$
After factorization of the $\sin$ function, it suffices by
\eqref{1.13} to show that
\begin{multline*}
\Bigl\|\,\int_0^\infty \cos(sP)F(s,\cd)\, ds\,
\Bigr\|_{\dot{H}^{\gamma-1}(\Omega)}+
\Bigl\|\,\int_0^\infty P^{-1}\sin(sP)F(s,\cd)\, ds\,
\Bigr\|_{\dot{H}^{\gamma}(\Omega)}
\\
\lesssim \|F\|_{L^{r'}_tY'(\R_+\times \Omega)}\,.
\end{multline*}
This, however, is the dual version
of \eqref{1.13} for $(Y,1-\gamma,r)$.
\qed

As a special case of \eqref{1.15} when the spaces $X$ and $Y$ are the 
standard Lebesgue spaces, we have the following

\begin{corr}\label{mixed}  Suppose that $n\ge 3$ and that Hypothesis B is 
valid.  Suppose that $q,\tilde q>2$, $r,\tilde r\ge 2$ and that
$$
\frac1q+\frac{n}r=\frac{n}2-\gamma=\frac1{\tilde q'}+\frac{n}{\tilde r'}-2
$$
and
$$
\frac2q+\frac{n-1}r\,, \, \frac2{\tilde q}+\frac{n-1}{\tilde r}\le \frac{n-1}2 .
$$
Then if the local Strichartz 
estimate \eqref{1.10} holds respectively for the triples
$\bigl(L^r(\Omega),\gamma,q\bigr)$ and
$\bigl(L^{\tilde r}(\Omega),1-\gamma,\tilde q\bigr)$,
it follows that when $u$ solves \eqref{1}
$$
\|u\|_{L^q_tL^r_x(\R_+\times\Omega)}
\lesssim \|f\|_{\dot H^\gamma_B(\Omega)}+\|g\|_{\dot H^{\gamma-1}_B(\Omega)}+
\|F\|_{L^{\tilde q'}_t\!L^{\tilde r'}_x(\R_+\times\Omega)}.
$$
These results also hold for $n=2$ under the above assumption, provided that 
$\gamma=1/2$.
\end{corr}

These estimates of course are the obstacle versions of the mixed-norm 
estimates for $\R^n$ and $\Delta_\g=\Delta$.  When $n\ge3$ (and \eqref{1.10} 
is valid) they include all the ones in the Keel-Tao theorem \cite{KT}, 
excluding the cases where either $q$ or $\tilde q$ is $2$.  For the 
Dirichlet-wave operator ($B=Id$) these results were proved in odd dimensions 
by Smith and Sogge \cite{SS} and then by Burq \cite{burq} and 
Metcalfe \cite{M} for even dimensions.  The Neumann case was not treated, 
but it follows from the same proof.  Unfortunately, the known techniques 
seem to only apply to the case of $\gamma=1/2$ when $n=2$, and Hypothesis B 
seems also to require $B=Id$ and $\partial\Omega\ne\emptyset$ in this case.  
The restriction that $\gamma=1/2$ when $n=2$ comes from the second part of 
(1.13), while for $n\ge3$ this is not an issue due to the fact that the 
Sobolev exponents $\gamma$ in Corollary 1.5 always satisfy 
$0\le \gamma\le 1$.  Also, at present, the knowledge of the local Strichartz 
estimates \eqref{1.10} when $X=L^r(\Omega)$ is limited.  
When $\Omega$ is the exterior of a geodesically convex obstacle, they were 
obtained by Smith and Sogge \cite{SS0}.  Recently, there has been work on 
proving local Strichartz estimates when $X=L^r(\Omega)$ for more general 
exterior domains (\cite{BLP}, \cite{BP}, \cite{BSS}, \cite{SS2}), but only 
partial results for a more restrictive range of exponents than the ones 
described in Corollary \ref{mixed} have been obtained.

%%%%%%%%%%%%%%%%%%%%%%%%%%%%%%%%%%%%%%%%%%%%%%%%%%%%%%%%

\newsection{Proof of Abstract Strichartz Estimates}

As mentioned before, we shall prove \eqref{1.13} by adapting the arguments 
from \cite{burq}, \cite{M} and \cite{SS}.  We shall assume that \eqref{1.2} 
is valid for $(\Omega, \Delta_\g)$ throughout.  A key step in the proof of
Theorem \ref{theorem1.3} will be to establish the following result that is 
implicit in \cite{burq}.

\begin{proposition}\label{prop1.5}  Let $w$ solve the inhomogeneous wave 
equation in Minkowski space
\begin{equation*}
\begin{cases}
(\partial_t^2-\Delta)w=F \quad \text{on} \; \R_+\times \R^n
\\
w|_{t=0}=\partial_tw|_{t=0}=0.
\end{cases}
\end{equation*}
Assume as above that \eqref{1.8} is valid whenever $v$ is a solution of 
the homogeneous wave equation \eqref{1.9}.
%, with, as above, $q$ and $\gamma$ satisfying \eqref{1.12}.  
Assume further that $q>2$ and $\gamma\ge -\tfrac{n-3}2$.
Then, if
$$
F(t,x)=0 \quad\mathrm{if}\quad |x|>2R\,,
$$
we have 
$$
\|w\|_{L^q_tX(\R_+\times \R^n)}\lesssim 
\|F\|_{L^2_t H^{\gamma-1}(\R_+\times \Rn)}.
$$
\end{proposition}

At the end of this section we shall show that when $n=2$ the assumption 
that $\gamma\ge1/2$ when $n=2$ is necessary even in the model case where 
$X=L^r(\R^n)$ with $2/q+1/r=1/2$ and $1/q+2/r=1-\gamma$.

To prove Proposition \ref{prop1.5}, we shall use our free space hypothesis 
\eqref{1.8} and 
the following result from \cite{SS}.

\begin{lemma}\label{lemma1.6}  Fix $\beta \in C^\infty_0(\Rn)$ and assume 
that $\gamma\le \tfrac{n-1}2$.  Then
$$\int_{-\infty}^\infty \Bigl\|\, \beta(\cd)\bigl(e^{it|D|}f\bigr)(t,\cd)\, 
\Bigr\|^2_{H^\gamma(\Rn)}\, dt \lesssim \|f\|_{\dot H^\gamma(\Rn)}^2,
$$
if $|D|=\sqrt{-\Delta}$.
\end{lemma}

As was shown in \cite{SS}, this lemma just follows from an application of 
Plancherel's theorem and the Schwarz inequality.  The assumption that 
$\gamma\le (n-1)/2$ is easily seen to be sharp.

To prove Proposition \ref{prop1.5}, we note that since we are assuming 
that $q>2$, by the Christ-Kiselev lemma \cite{CK}, it suffices to show that
\begin{equation}\label{11}
\Bigl\| \, \int_{0}^\infty e^{i(t-s)|D|} 
|D|^{-1} \beta(\cd)G(s,\cd)\, ds \, \Bigr\|_{L^q_tX(\R_+\times\R^n)}
\lesssim \|G\|_{L^2_tH^{\gamma-1}(\Rn)},
\end{equation}
assuming that $\beta\in C^\infty_0(\Rn)$.  If we apply \eqref{1.8}, 
we conclude that the left side of this inequality is majorized by
$$
\Bigl\| \, \int_{0}^\infty e^{-is|D|} 
|D|^{-1+\gamma} \beta(\cd)G(s,\cd)\, ds \, 
\Bigr\|_{L^2(\Rn)}.
$$
Since $\| (1-\Delta)^{(\gamma-1)/2}G(s,\cd)\|_2=\|G(s,\cd)\|_{H^{\gamma-1}}$, 
it suffices to see that
$$
\Bigl\|\, \int_{0}^\infty e^{-is|D|}
|D|^{-1+\gamma}\beta(\cd)(1-\Delta)^{(1-\gamma)/2}H(s,\cd)\, ds \, 
\Bigr\|_{L^2(\R^n)}\lesssim \|H\|_{L^2(\R_+\times\R^n)}.
$$
By duality, this is equivalent to the statement that
\begin{equation}\label{12}
\bigl\|\, (1-\Delta)^{(1-\gamma)/2}\beta(\cd)e^{is|D|}
|D|^{-1+\gamma}h\bigr\|_{L^2(\R_+\times\R^n)}
\lesssim \|h\|_{L^2(\R^n)}.
\end{equation}
Since we are assuming that $\gamma\ge -\tfrac{n-3}2$, 
we have that $1-\gamma \le \tfrac{n-1}2$.
Therefore, \eqref{12} follows from Lemma \ref{lemma1.6}, completing 
the proof of Proposition \ref{prop1.5}.
\qed

To prove Theorem \ref{theorem1.3} we also need a similar result for 
solutions of the wave equation \eqref{1} for $(\Omega, B, \Delta_\g)$.

\begin{proposition}\label{prop1.7}  Let $u$ solve \eqref{1} and assume that
\begin{equation}\label{13}
f(x)=g(x)=F(t,x)=0,\quad \text{when } \, |x|>2R.
\end{equation}
Then if $(X,\gamma,q)$ is an admissible triple with 
%$(\gamma,q)$ 
%as in \eqref{1.12} we have
$q>2$ and $\gamma\ge -\tfrac{n-3}2$ we have
\begin{equation}\label{14}
\|u\|_{L^q_tX(\R_+\times\Omega)}\lesssim 
\|f\|_{ H^\gamma_B}+\|g\|_{H^{\gamma-1}_B}+\|F\|_{L^2_tH^{\gamma-1}_B}.
\end{equation}
\end{proposition}

The key ingredients in the proof are Proposition \ref{prop1.5} and the 
following variant of \eqref{1.2}, which holds for all $\gamma\in \R$,
provided \eqref{13} holds, and
$\beta\in C_c^\infty(\Rn)$ equals 1 on a neighborhood of 
$\R^n\backslash\Omega$:
\begin{multline}\label{22}
\|\beta u\|_{L^\infty_tH^\gamma_B}+
\|\beta\partial_t u\|_{L^\infty_tH^{\gamma-1}_B}+
\|\beta u\|_{L^2_tH^\gamma_B}+\|\beta\partial_t u\|_{L^2_tH^{\gamma-1}_B}
\\
\lesssim \|f\|_{H^\gamma_B}+
\|g\|_{H^{\gamma-1}_B}+\|F\|_{L^2_tH^{\gamma-1}_B}.
\end{multline}
The $L^2_t$ estimates in \eqref{22}
on $u$ follow from \eqref{1.2} and elliptic regularity arguments for 
$\gamma\in\Z$, and by interpolation for the remaining $\gamma\in \R$.
The $L^\infty_t$ estimates then follow from energy estimates, duality, and
elliptic regularity.

To prove \eqref{14}, let us fix $\beta\in C^\infty_0(\Rn)$ satisfying 
$\beta(x)=1$, $|x|\le 3R$ and write
$$u=v+w, \quad \text{where } \, v=\beta u, \, \, w=(1-\beta)u.$$
Then $w$ solves the free wave equation
\begin{equation*}
\begin{cases}
(\partial_t^2-\Delta)w=[\beta,\Delta]u
\\
w|_{t=0}=\partial_tw|_{t=0}=0.
\end{cases}
\end{equation*}
An application of Proposition \ref{prop1.5} shows that $\|w\|_{L^q_tX}$ is 
dominated by 
$\|\rho u\|_{L^2_tH^\gamma_B}$ if $\rho$ equals one on the support of $\beta$.
Therefore, by \eqref{22},
%Lemma \ref{lemma1.1}, 
$\|w\|_{L^q_tX}$ is dominated by the right side of \eqref{14}.

As a result, we are left with showing that if $v=\beta u$ then
\begin{equation}\label{15}
\|v\|_{L^q_tX(\R_+\times \Omega)}\lesssim 
\|f\|_{H^\gamma_B}+\|g\|_{ H^{\gamma-1}_B}+\|F\|_{L^2_tH^{\gamma-1}_B},
\end{equation}
assuming, as above, that \eqref{13} holds.  To do this, fix 
$\varphi \in C^\infty_0((-1,1))$ satisfying 
\linebreak
$\sum_{j=-\infty}^\infty \varphi(t-j)=1$.  For a given 
$j\in {\mathbb N}$, let $v_j=\varphi(t-j)v$.  Then $v_j$ solves
\begin{equation*}
\begin{cases}
(\partial_t^2-\Delta_\g)v_j = -\varphi(t-j)[\Delta,\beta]u + 
[\partial^2_t,\varphi(t-j)]\beta u
+\varphi(t-j)F
\\
Bv_j(t,x)=0,\quad x\in \partial\Omega
\\
v_j(0,\cd)=\partial_tv_j(0,\cd)=0,
\end{cases}
\end{equation*}
while $v_0=v-\sum_{j=1}^\infty v_j$ solves
\begin{equation*}
\begin{cases}
(\partial_t^2-\Delta_\g)v_0=-\tilde \varphi [\Delta,\beta]u + 
[\partial_t^2,\tilde \varphi]\beta u + \tilde \varphi F
\\ 
Bv_0(t,x)=0,\quad x\in \partial\Omega
\\
v_0|_{t=0}=f,\, \, \partial_t v_0|_{t=0}=g,
\end{cases}
\end{equation*}
if $\tilde \varphi =1-\sum_{j=1}^\infty \varphi(t-j)$ if $t\ge 0$ and $0$ 
otherwise.
If we then let $G_j=(\partial_t^2-\Delta_\g)v_j$ be the forcing term for 
$v_j$, $j=0,1,2,\dots$, then, by 
\eqref{22},
%Lemma \ref{lemma1.1},
 we have that
$$\sum_{j=0}^\infty \|G_j\|^2_{L^2_tH^{\gamma-1}_B(\R_+\times \Omega)}
\lesssim \|f\|^2_{H^{\gamma}_B}+\|g\|^2_{H^{\gamma-1}_B}+
\|F\|_{L^2_t H^{\gamma-1}_B}^2.$$
By the local Strichartz estimates \eqref{1.10} and Duhamel, we get for 
$j=1,2,\dots$
$$\|v_j\|_{L^q_tX(\R_+\times \Omega)}\lesssim
\int_0^\infty \|G_j(s,\cd)\|_{ H^{\gamma-1}_B}\, ds
\lesssim \|G_j\|_{L^2_t H^{\gamma-1}_B},$$
using Schwarz's inequality and the support properties of the $G_j$ in the 
last step.  Similarly, 
$$
\|v_0\|_{L^q_tX(\R_+\times \Omega)}\lesssim 
\|f\|_{H^\gamma_B}+\|g\|_{H^{\gamma-1}_B}+\|G_0\|_{L^2_tH^{\gamma-1}_B}\,.
$$
Since $q>2$, we have
$$
\|v\|^2_{L^q_tX(\R_+\times \Omega)}\lesssim 
\sum_{j=0}^\infty \|v_j\|^2_{L^q_tX(\R_+\times \Omega)}
$$
and so we get
$$
\|v\|^2_{L^q_tX}\lesssim 
\|f\|^2_{H^\gamma_B}+\|g\|^2_{H^{\gamma-1}_B}+\|F\|_{L^2_tH^{\gamma-1}}^2
$$
as desired, which finishes the proof of Proposition \ref{prop1.7}.
\qed

\medskip

\noindent{\bf End of Proof of Theorem \ref{theorem1.3}:}
Recall that we are assuming that $(\partial_t^2-\Delta_\g)u=0$.  
By Proposition \ref{prop1.7} we may also assume that the initial data 
for $u$ vanishes when $|x|<3R/2$.  We then fix $\beta\in C^\infty_0(\Rn)$ 
satisfying $\beta(x)=1$, $|x|\le R$ and $\beta(x)=0$, $|x|>3R/2$ and write
$$
u=u_0-v = (1-\beta)u_0 + ( \beta u_0-v)\,,
$$
where $u_0$ solves the Cauchy problem for the Minkowski space wave equation 
with initial data defined to be $(f,g)$ if $x\in \Omega$ and $0$ otherwise.  
By the free estimate \eqref{1.8} and \eqref{local}, 
we can restrict our attention to 
$\tilde u= \beta u_0-v$.
But
$$
(\partial^2_t-\Delta_\g)\tilde u = -[\Delta, \beta]u_0\equiv G
$$
is supported in $R<|x|<2R$, and satisfies
\begin{equation}\label{Gdecay}
\int_{0}^\infty \|G(t,\cd)\|^2_{H^{\gamma-1}_B}\, dt 
\lesssim \|f\|_{\dot H^\gamma_B}^2+\|g\|_{\dot H^{\gamma-1}_B}^2
\end{equation}
by Lemma \ref{lemma1.6} and the fact that $G$ vanishes on a neighborhood
of $\partial\Omega$.
Note also that $\tilde u$ has vanishing initial data.
Therefore, since Proposition \ref{prop1.7} tells us that 
$\|\tilde u\|_{L^q_tX(\R_+\times\R^n)}^2$ 
is dominated by the left side of \eqref{Gdecay}, 
the proof is complete.
\qed

For future reference, we note that the preceeding steps
establish the following generalization of \eqref{22},
assuming that $\gamma\in[-\frac{n-3}2,\frac{n-1}2]$, and that
$F(x)=0$ for $|x|>R$:

\begin{multline}\label{22'}
\|u\|_{L^\infty_t\dot{H}^\gamma_B}+
\|\partial_t u\|_{L^\infty_t\dot{H}^{\gamma-1}_B}+
\|\beta u\|_{L^2_tH^\gamma_B}+\|\beta\partial_t u\|_{L^2_tH^{\gamma-1}_B}
\\
\lesssim \|f\|_{\dot{H}^\gamma_B}+
\|g\|_{\dot{H}^{\gamma-1}_B}+\|F\|_{L^2_tH^{\gamma-1}_B}.
\end{multline}
In particular, $f$ and $g$ have no support restrictions.
To see that \eqref{22'} holds, first consider bounding the terms
$\|\beta\partial_t^j u\|_{L^2_tH_B^{\gamma-j}}$ for $j=0,1$. 
For these terms, it suffices by \eqref{22} to consider $F=0$ 
and $f,g=0$ near $\partial\Omega$. Decomposing $u=(1-\beta)u_0+\tilde u$
as above, we may use \eqref{22} and \eqref{Gdecay} to deduce 
the $L^2_t$ bounds in \eqref{22'} for $u$.
These bounds now yield
$$
\|(\partial_t^2-\Delta_\g)(1-\beta)u\|_{L^2_tH^{\gamma-1}_B}+
\|(\partial_t^2-\Delta_\g)\beta u\|_{L^2_tH^{\gamma-1}_B}
\lesssim
\|f\|_{\dot{H}^\gamma_B}+
\|g\|_{\dot{H}^{\gamma-1}_B}+\|F\|_{L^2_tH^{\gamma-1}_B}\,.
$$
The $L^\infty_t$ bounds on $\beta u$ now follow from \eqref{22}.
Finally, $(1-\beta)u$ satisfies the Minkowski wave equation on
$\R\times\R^n$, with initial data in $\dot H^\gamma\times\dot H^{\gamma-1}$,
and driving force $\tilde F \in L^2_t\dot H^{\gamma-1}$ which vanishes 
for $|x|\ge R$.
The contribution to $u$ from its initial data satisfies the $L^\infty_t$
bounds as a result of homogeneous Sobolev bounds for the Minkowski wave group.
The contribution from $\tilde F$ 
is bounded using Lemma \ref{lemma1.6} and duality.

\medskip

Let us conclude this section by showing that when $n=2$ the restriction in 
Proposition~\ref{prop1.5} that $\gamma\ge 1/2$ is necessary in the case where 
$X=L^r(\R^2)$.  In this case, by the standard mixed-norm Strichartz estimates 
(see e.g. \cite{KT}), the hypotheses of the Proposition are satisfied when 
$0\le \gamma<3/4$, $1/q+2/r=1-\gamma$ and $2/q+1/r \le 1/2$.  
%For later reference, we note that the last condition 
%just says that $q=4r/(r-2)$ in what follows.

Since the hypotheses are satisfied, if the Proposition were valid for a 
given $\gamma$ and $X=L^r(\R^2)$ as above, then the 
$L^q_tL^r_x(\R_+\times \R^2)$ norm of 
$$WF(t,x)=\int_{-\infty}^t \int_{\R^2} e^{ix\cdot \xi}
\frac{\sin (t-s)|\xi|}{|\xi|} \Hat F(s,\xi)\, d\xi ds
$$
would have to be bounded by the $L^2_tH^{\gamma-1}$ norm of $F$ if 
$F(t,x)=0$ when $|x|>1$.
We shall take $F$ to be a product $h_T(s)\beta(x)$ where 
$\beta\in C^\infty(\R^2)$ vanishes for 
 $|x|>1$ but satisfies $\hat \beta(0)=1$, while $h_T$ is an odd 
function supported in $[-T,T]$.  For this choice of $F$ we have
$$WF(t,x)=-i \int_{\R^2}e^{ix\cdot\xi}\cos(t|\xi|) \Hat h_T(|\xi|)
\Hat \beta(\xi)\, d\xi/|\xi|, \quad \text{if } t>T.$$
Fix a nonzero function $\rho\in C^\infty(\R)$ supported in $(1/2,1)$.  
If we take $h_T$ to be the odd function which equals $T^{-1/2}\rho(s/T)$ 
for positive $s$, then since $h_T$ has a non-zero $L^2$ norm which is 
independent of $T$, if Proposition~\ref{prop1.5} were valid for an 
$L^q_tL^r_x$ space as above, then it would follow that
\begin{multline*}WF(t,x)=
-iT^{1/2}\int_{\R^2} e^{ix\cdot \xi} \cos(t|\xi|) \Hat h_1(T|\xi|)
\Hat \beta(\xi)\, d\xi/|\xi|
\\ 
=-iT^{-1/2}
\int_{\R^2}e^{i\frac{x}{T}\cdot \xi} \cos(\tfrac{t}{T}  |\xi|) \, 
\Hat h_1(|\xi|) \, \Hat \beta(\xi/T) \, d\xi/|\xi|
\end{multline*}
would belong to $L^q_tL^r_x([T,\infty)\times \R^2)$ with a bound 
independent of $T$.  An easy calculation shows that this norm equals
$$T^{-1/2+1/q+2/r}\Bigl\|\int_{\R^2}e^{ix\cdot\xi}\cos(t|\xi|) 
\Hat h_1(|\xi|)\Hat \beta(\xi/T)\frac{d\xi}{|\xi|}\, 
\Bigr\|_{L^q_tL^r_x([1,\infty)\times \R^2)}.$$
Since our assumption that $\Hat \beta (0)=1$ implies that the last 
factor on the right tends to a positive constant, we conclude that if 
the conclusion of Proposition~\ref{prop1.5} were valid for $X=L^r(\R^2)$, 
then we would need that
%$$-\tfrac12+\tfrac1q+\tfrac2r=-\tfrac12+\tfrac{r-2}{4r}+\tfrac2r\le 0,$$
%which is equivalent to $r\ge6$.  Since $1-\gamma=1/q+2/r=(r-2)/4r+2/r$, 
%this is equivalent to $\gamma\ge 1/2$, as claimed.  
$$\frac12 \ge \frac1q+\frac2r=1-\gamma.$$
This means that when $n=2$, the assumption that  $\gamma\ge 1/2$ in 
Proposition~\ref{prop1.5} is necessary.

%%%%%%%%%%%%%%%%%%%%%%%%%%%%%%%%%%%%%%%%%%%%%%%%%%%%%%%%

\newsection{The Strauss conjecture for nontrapping obstacles when $n=3,4$}

Let us start the proof of Theorem \ref{exist} by going over the Minkowski 
space results that will be used.  These will form the assumption \eqref{1.8} 
of Theorem \ref{theorem1.3}.

\begin{lemma}
  Let $u$ solve the Minkowski wave equation
   $$(\partial_t^2-\Delta)u=F,\quad (t,x)\in \R\times\R^n$$
   $$u(0,\cd)=f,\quad \partial_tu(0,\cd)=g.$$
Then, for $2\le p\le \infty$, and $\gamma$ satisfying
\begin{equation}\label{gamma_restrictions}
  \frac{1}{2}-\frac{1}{p}<\gamma<\frac{n}{2}-\frac{1}{p},\quad\text{
    and }\quad \frac{1}{2}<1-\gamma<\frac{n}{2},
\end{equation}
  we have the following estimate
  \begin{multline}\label{Minkowski_weighted_Strichartz}
    \Bigl\||x|^{\frac{n}{2}-\frac{n+1}{p}-\gamma}u
\Bigr\|_{L^p_tL^p_rL^2_\omega(\R_+\times\R^n)}\lesssim
\|f\|_{\dot{H}^\gamma(\R^n)}+\|g\|_{\dot{H}^{\gamma-1}(\R^n)} \\+
\Bigl\||x|^{-\frac{n}{2}+1-\gamma}F\Bigr\|_{L^1_tL^1_rL^2_\omega(\R_+\times \R^n)}.
  \end{multline}
\end{lemma}

Here, and in what follows, we are using the mixed-norm notation with respect 
to the volume element
$$
\|h\|_{L^q_rL^p_\omega}=\Bigl(\, \int_0^\infty \, 
\Bigl(\, \int_{S^{n-1}}|h(r\omega)|^p\, d\sigma(\omega)\, \Bigr)^{q/p}
\, r^{n-1}dr\, \Bigr)^{1/q}
$$
for finite exponents and
$$
\|h\|_{L^\infty_rL^p_\omega}=\sup_{r>0}
\Bigl(\, \int_{S^{n-1}}|h(r\omega)|^p\, d\sigma(\omega)\, \Bigr)^{1/p}.
$$

We first note that, by the trace lemma for the unit sphere 
and scaling, we have
\begin{equation}\label{3.1}
\sup_{r>0}r^{\frac{n}2-s}\Bigl(\, \int_{S^{n-1}}|v(r\omega)|^2 \, 
d\sigma(\omega)\, \Bigr)^{1/2}\lesssim \|v\|_{\dot H^s(\R^n)}\,, 
\quad \frac12<s<\frac{n}2\,,
\end{equation}
where $d\sigma$ denotes the unit measure on $S^{n-1}$.  Consequently,
$$\sup_{r>0}r^{\frac{n}2-s}\Bigl(\int_{S^{n-1}}\bigl| \, 
\bigl(e^{it|D|}\varphi \bigr)(r\omega)\, \bigr|^2\, d\sigma(\omega)\, 
\Bigr)^{1/2}
\lesssim \|\varphi \|_{\dot H^{s}(\R^n)}\,, \quad \frac12<s<\frac{n}2\,,$$
which is equivalent to 
\begin{equation}\label{3.2}
\|\, |x|^{-\alpha} e^{it|D|}\varphi\|_{L^\infty_rL^2_\omega}
\lesssim \|\varphi\|_{\dot H^{\frac{n}2+\alpha}(\R^n)}, \quad  
-\frac{n-1}2<\alpha<0.
\end{equation}

Note that by applying \eqref{3.1} to the Fourier transform of $v$, 
we see that it is equivalent to the uniform bounds
$$
\Bigl(\, \int_{S^{n-1}}|\Hat v(\lambda\omega)|^2\, d\sigma(\omega)\, 
\Bigr)^{1/2}\lesssim
\lambda^{-\frac{n}2+s}\|\, |x|^s v\|_{L^2(\Rn)}\,, 
\quad \lambda>0\,, \quad
\frac12<s<\frac{n}2\,,
$$
which by duality is equivalent to 
\begin{equation}\label{3.3}
\Bigl\|\, |x|^{-s}\int_{S^{n-1}}h(\omega) e^{i\lambda x\cdot \omega}\, 
d\sigma(\omega) \,  \Bigr\|_{L^2_x(\R^n)}\lesssim 
\lambda^{s-\frac{n}2}\|h\|_{L^2_\omega(S^{n-1})}\,,
\end{equation}
for $\lambda>0$ and fixed $1/2<s<n/2$.  Using this estimate we can obtain
\begin{equation}\label{3.4}
\bigl\|\, |x|^{-s} e^{it|D|}\varphi \|_{L^2(\R_+ \times \R^n)} \lesssim 
\|\, |D|^{s-\frac12}\varphi\|_{L^2(\R^n)}, \quad \frac12<s<\frac{n}2\,,
\end{equation}
for, by after Plancherel's theorem with respect to the $t$-variable, we find 
that the square of the left side of \eqref{3.4} equals 
\begin{multline*}
(2\pi)^{-1} \int_0^\infty \int_{\R^n}\Bigl| \, |x|^{-s}\int_{S^{n-1}} 
e^{ix\cdot \rho \omega}\rho^{n-1}\Hat \varphi(\rho\omega)\, d\sigma(\omega)\, 
\Bigr|^2\, dx \, d\rho
\\
\lesssim \int_0^\infty \int_{S^{n-1}}\rho^{2(n-1)}|
\Hat \varphi(\rho\omega)|^2 \, \rho^{2s-n}\, d\sigma(\omega)d\rho = 
\|\, |D|^{s-\frac12}\varphi\|_{L^2(\R^n)}^2,
\end{multline*}
using \eqref{3.3} in the first step.

If we interpolate between \eqref{3.2} and \eqref{3.4} we conclude that,
for $2\le q\le \infty$,
\begin{equation}\label{3.5}
\Bigl\|\, |x|^{\frac{n}2-\frac{n+1}q-\gamma} e^{it|D|}\varphi \, 
\Bigr\|_{L^q_tL^q_rL^2_\omega(\R_+\times \Rn)}
\lesssim \|\varphi\|_{\dot H^\gamma(\Rn)}, \quad 
\frac12-\frac1q<\gamma<\frac{n}2-\frac1q.
\end{equation}
This estimate in turn implies that if $v$ solves the Cauchy problem 
$(\partial_t^2-\Delta)v=0$ in $\R_+\times \Rn$ then
\begin{multline}\label{3.6}
\Bigl\|\, |x|^{\frac{n}2-\frac{n+1}q-\gamma}v \, 
\Bigr\|_{L^q_tL^q_rL^2_\omega(\R_+\times \Rn)}
\\
\lesssim \|v(0,\cd)\|_{\dot H^\gamma(\Rn)}+
\|\partial_t v(0,\cd)\|_{\dot H^{\gamma-1}(\Rn)}, \, \, \, 
\frac12-\frac1q <\gamma<\frac{n}2-\frac1q\,.
\end{multline}
The estimate dual to \eqref{3.1} is
\begin{equation}\label{dual}
\|\varphi\|_{\dot H^{\gamma-1}}\le 
\bigl\|\,|x|^{-\frac n2+1-\gamma}\varphi\bigr\|_{L^1_rL^2_\omega}\,.
\end{equation}
By the Duhamel formula and \eqref{3.6}-\eqref{dual}, we then have
\begin{multline}\label{3.7}
\Bigl\|\, |x|^{\frac{n}2-\frac{n+1}p-\gamma} u\, 
\Bigr\|_{L^p_tL^p_rL^2_\omega(\R_+\times \Rn)}
\lesssim \|u(0,\cd)\|_{\dot H^\gamma(\Rn)}+
\|\partial_t u(0,\cd)\|_{\dot H^{\gamma-1}(\Rn)}
\\
+\Bigl\|\, |x|^{-\frac{n}2+1-\gamma} (\partial_t^2-\Delta)u\, 
\Bigr\|_{L^1_tL^1_rL^2_\omega(\R_+\times\Rn)},
\end{multline}
provided that $\gamma$ and $1-\gamma$ satisfy the condition in \eqref{3.6} 
for $q$ equal to $p$ and $\infty$, respectively, i.e., 
\eqref{gamma_restrictions}.\qed
%\begin{equation}\label{3.8}
%\frac12-\frac1p<\gamma<\frac{n}2-\frac1p \quad \text{and} \quad 
%\frac12<1-\gamma<\frac{n}2.
%\end{equation}\qed

A calculation shows that if 
\begin{equation}\label{3.9} \gamma=\frac{n}2-\frac2{p-1},
\end{equation}
and
$$p_c<p< (n+3)/(n-1)$$
then \eqref{gamma_restrictions} holds: $p>p_c$ is needed for the first part, 
and 
$p<(n+3)/(n-1)$ for the second.  Additionally, as far as the powers of 
$|x|$ go in \eqref{3.7}, we have
\begin{equation}\label{3.10}
p\Bigl(\, \frac{n}2-\frac{n+1}p-\gamma\, \Bigr)=
p\Bigl(\, \frac{(n+1)-(n-1)p}{p(p-1)}\, \Bigr)
= -\frac{n}2+1-\gamma\,, \quad \text{if } \, \, 
\gamma=\frac{n}2-\frac2{p-1}.
\end{equation}
As a result, by the arguments to follow, \eqref{Minkowski_weighted_Strichartz} 
is strong enough to show 
that for the non-obstacle, Minkowski space case, i.e.~ $\Omega=\Rn$, 
$\Delta_\g=\Delta$, if $2\le n\le 4$ then for $p_c<p<(n+3)/(n-1)$, the 
equation \eqref{1.3} has a global solution for small data as described in 
Theorem \ref{theorem1.3}.

To prove the obstacle version of this result for $n=3$ and $4$ we shall use 
a slightly weaker inequality for which it will be easy to show that we have 
the corresponding local Strichartz estimates \eqref{1.10} for 
$(\Omega,\Delta_\g)$.  To this end, if $R$ is chosen so that $\partial\Omega$ 
is contained in $|x|<R$ and $\Delta=\Delta_\g$ for $|x|\ge R$ then we define 
$X=X_{\gamma,q}(\Rn)$ to be the space with norm defined by
\begin{equation}\label{3.11}
\|h\|_{X_{\gamma,q}}=\|h\|_{L^{s_\gamma}(|x|<2R)}\, + 
\, \bigl\| \, |x|^{\frac{n}2-\frac{n+1}q-\gamma}h\|_{L^q_rL^2_\omega(|x|>2R)}
, \quad \text{if }\, \, n\bigl(\tfrac12-\tfrac1{s_\gamma})=\gamma.
\end{equation}
We then prove the following obstacle variant of 
\eqref{Minkowski_weighted_Strichartz}.

\begin{lemma}\label{weighted_Strichartz}
For solutions of \eqref{1} if $n\ge 3$ and $p>2$:
\begin{multline}\label{3.13}
\Bigl\|\, |x|^{\frac{n}2-\frac{n+1}p-\gamma}u\, 
\Bigr\|_{L^p_tL^p_rL^2_\omega(\R_+\times
\{|x|>2R\})}\, + \, 
\|u\|_{L^p_tL^{s_\gamma}_x(\R_+\times \{x\in \Omega: \, |x|<2R\})}
\\
\lesssim \|f\|_{\dot H^\gamma_B}+
\|g\|_{\dot H^{\gamma-1}_B}
+\Bigl\|\, |x|^{-\frac{n}2+1-\gamma}F
\Bigr\|_{L^1_tL^1_rL^2_\omega(\R_+\times\{|x|>2R\})}
\\
+ \|F\|_{L^1_tL^{s_{1-\gamma}'}_x(\R_+\times \{x\in \Omega: |x|<2R\})}
\end{multline}
provided that \eqref{gamma_restrictions} holds.
\end{lemma}

By \eqref{3.6} and Lemma \ref{lemma1.6} we have that the assumption 
\eqref{1.8} of Theorem \ref{theorem1.3} is valid if $1/2-1/q<\gamma<n/2-1/q$ 
and $2\le q\le \infty$, i.e.
\begin{multline}\label{3.12}
\|v\|_{L^q_t X_{\gamma,q}(\R_+\times \Rn)}\lesssim 
\|v(0,\cd)\|_{\dot H^\gamma(\Rn)}
+\|\partial_tv(0,\cd)\|_{\dot H^{\gamma-1}(\Rn)},
\\
\text{if } \, \, (\partial_t^2-\Delta)v=0 \, \, \, \text{in } \, \, 
\R_+\times \Rn,
\end{multline}
under the additional assumption that $\gamma\le (n-1)/2$ (which is the case 
for \eqref{3.9}).  Indeed the contribution of the second part of the norm in 
\eqref{3.11} is controlled by \eqref{3.6}.  To handle the contribution of the 
first term in the right side of \eqref{3.11} we note that if 
$\beta\in C^\infty_0(\Rn)$ equals one when $|x|\le 3R$ then Sobolev estimates 
yield
$$\|v(t,\cd)\|_{L^{s_\gamma}(|x|<2R)}\lesssim 
\|\beta(\cd)v(t,\cd)\|_{\dot H^\gamma(\Rn)}.$$
Thus, $\|v\|_{L^q_tL^{s_\gamma}(\R_+\times \{|x|<R\})}$ is controlled by 
the right side of \eqref{3.12} for $q=2$, by Lemma \ref{lemma1.6}.  Since 
this is also the case for $q=\infty$ by energy estimates, by interpolation 
we conclude that we can control the contribution of the first term in the 
right side of \eqref{3.11} to \eqref{3.12}, which finishes the proof
of \eqref{3.12}.

Since the dual norm of $\|\, |x|^{\alpha}h\|_p$ is 
$\|\, |x|^{-\alpha}h\|_{p'}$, by Corollary 1.4, we would get \eqref{3.13} 
from \eqref{3.12} and Hypothesis B if we could show that for $q>2$
$$\|u\|_{L^q_t X_{\gamma,q}([0,1]\times \Omega)}
\lesssim \|f\|_{\dot H^\gamma_B}+
\|g\|_{\dot H^{\gamma-1}_B},
$$
whenever $u$ solves \eqref{1} with $F\equiv 0$, and, as above, 
$1/2-1/q<\gamma<n/2-1/q$.  By the finite propagation speed of the wave 
equation, it is clear that the contribution of the second term in the right 
side of \eqref{3.11} will enjoy this estimate.  As before, the first term 
satisfies it because of Sobolev estimates.  This completes the 
proof of \eqref{3.13}. \qed

Let us also observe a related estimate
\begin{multline}\label{a.14}%\label{a.17}
\|u\|_{L^\infty_t\dot H^\gamma_B(\R_+\times\Omega)}+ 
\|\partial_tu\|_{L^\infty_t\dot H^{\gamma-1}_B(\R_+\times\Omega)} + 
\|u\|_{L^\infty_tL^{s_\gamma}_x(\R_+\times\Omega)}+
\|\beta u\|_{L^2_tH^\gamma_B(\R_+\times \Omega)}
\\
\lesssim \|f\|_{\dot H^\gamma_B}+
\|g\|_{\dot H^{\gamma-1}_B}
+\Bigl\|\, |x|^{-\frac{n}2+1-\gamma}F\Bigr\|_{L^1_tL^1_rL^2_\omega(\R_+\times
\{|x|>2R\})}
\\
+ \|F\|_{L^1_tL^{s_{1-\gamma}'}_x(\R_+\times \{x\in \Omega: |x|<2R\})},
\end{multline}
assuming that \eqref{gamma_restrictions} holds.
Indeed, this is a direct consequence of \eqref{22'} and the Duhamel
formula, together with the inclusion 
$\dot H^\gamma_B(\Omega)\hookrightarrow L^{s_\gamma}(\Omega)$, and the
following consequence of \eqref{dual}, and the dual estimate to
Sobolev embedding 
$\dot H^{1-\gamma}_B(\Omega)\hookrightarrow L^{s_{1-\gamma}}(\Omega)$,
\begin{equation}\label{a.14'}
\|g\|_{\dot{H}^{\gamma-1}_B}\lesssim
\|\, |x|^{-\frac{n}2+1-\gamma}g\|_{L^1_rL^2_\omega(|x|>2R)}
+\|g\|_{L^{s'_{1-\gamma}}(x\in \Omega: |x|<2R)}\,.
\end{equation}

To prove Theorem \ref{exist} we shall require a variation of the last two 
estimates involving the vector fields
$$\{\Gamma\}=\{\partial_t,Z\}$$
where, as before, $\{Z\}$ are the vector fields 
$\{\partial_i,\, x_j\partial_k-x_k\partial_j: 
\, 1\le i\le n, 1\le j<k\le n\}$.  Note that all the $\{\Gamma\}$ commute 
with $\square_g=\partial_t^2-\Delta_\g$ when $|x|>R$ because 
$\partial\Omega \subset \{x: \, |x|<R\}$ and $\Delta = \Delta_\g$ for $|x|>R$.

The main estimate we require is the following.
\begin{lemma}
With $p$ and $\gamma$ as in Lemma \ref{weighted_Strichartz}, $u$
solving
\eqref{1} with $n\ge 3$, and $(f,g,F)$ satisfying 
$H^{2}_B\times H^{1}_B\times H^{1}_B$ 
boundary conditions, then
\begin{align}\label{a.15}
&\sum_{|\alpha|\le 2}
\Bigl( \|\, |x|^{\frac{n}2-\frac{n+1}p-\gamma} \Gamma^\alpha u
\|_{L^p_tL^p_rL^2_\omega(\R_+\times \{|x|>2R\})}
+\|\Gamma^\alpha u\|_{L^p_tL^{s_\gamma}_x(\R_+\times \{x\in \Omega: |x|<2R\})}
\Bigr)
\\
&\lesssim \sum_{|\alpha|\le 2}
\Bigl( \|Z^\alpha f\|_{\dot H^\gamma_B}
+\|Z^\alpha g\|_{\dot H^{\gamma-1}_B}\Bigr)
\notag
\\
&+\sum_{|\alpha|\le 2}
\Bigl(\|\, |x|^{-\frac{n}2+1-\gamma}
\Gamma^\alpha F\|_{L^1_tL^1_rL^2_\omega(\R_+\times \{|x|>2R\})}+
\|\Gamma^\alpha F\|_{L^1_tL_x^{s'_{1-\gamma}}
(\R_+\times \{x\in \Omega: |x|<2R\})}\Bigr).\notag
%\\
%&\qquad+\sum_{|\alpha|\le 1}\|\Gamma^\alpha F\|_{L^p_tL_x^{s'_{1-\gamma}}
%(\R_+\times \{x\in \Omega: \, |x|<2R\})}.\notag
\end{align}
\end{lemma}

The boundary conditions on $(f,g,F)$ imply that $\partial_t^j u$
is locally in $H^{2+\gamma-j}(\Omega)$, $j=0,1,2,$ which will be
implicitly used in elliptic regularity arguments. 
We will also use the fact that the Cauchy data for $\Gamma^\alpha u$ is
bounded in $\dot H^\gamma_B\times \dot H^{\gamma-1}_B$ by the right hand side
of \eqref{a.15} for $|\alpha|\le 2$. This is clear if $\Gamma^\alpha$ is 
replaced
by $Z^\alpha$. On the other hand, the Cauchy data for $\partial_t u$ is
$(g,\Delta_\g f+F(0,\cd))$. We may control
$$
\sum_{|\alpha|\le 1}\Bigl(\,\|Z^\alpha g\|_{\dot H^\gamma_B}+
\|Z^\alpha \Delta_\g f\|_{\dot H^{\gamma-1}_B}\Bigr)\le
\sum_{|\alpha|\le 2}
\Bigl(\,\|Z^\alpha f\|_{\dot H^\gamma_B}
+\|Z^\alpha g\|_{\dot H^{\gamma-1}_B}\Bigr)\,.
$$
Recall that $\gamma\in (0,\frac 12)$, so that 
$\dot H^\gamma_B(\Omega)=\dot H^\gamma(\overline{\Omega})$.
To control the term $F(0,\cd)$, we recall that $\Gamma=\{\partial_t,Z\}$,
and use the bound
\begin{equation}\label{a.15'}
\sum_{|\alpha|\le 1}\|\Gamma^\alpha F
\|_{L^\infty_t\dot H^{\gamma-1}_B(\R_+\times\Omega)}
\le
\sum_{|\alpha|\le 2}\|\Gamma^\alpha F
\|_{L^1_t\dot H^{\gamma-1}_B(\R_+\times\Omega)}
\end{equation}
which by \eqref{a.14'} is seen to be dominated by the 
right hand side of \eqref{a.15}.
Similar considerations apply to the Cauchy data for $\partial_t^2 u$.

Let us now give the argument for \eqref{a.15}.  We first fix 
$\beta_0\in C^\infty_0$ satisfying $\beta_0=1$ for $|x|\le R$ and 
$\supp \beta_0\subset \{|x|<2R\}$.  Then the first step in the proof of 
\eqref{a.15} will be to show that
\begin{align}\label{a.18}
&\sum_{|\alpha|\le 2}\Bigl( \|\, |x|^{\frac{n}2-\frac{n+1}p-\gamma}(1-\beta_0) 
\Gamma^\alpha u\|_{L^p_tL^p_rL^2_\omega(\R_+\times \{|x|>2R\})}
+\| (1-\beta_0)\Gamma^\alpha u
\|_{L^p_tL^{s_\gamma}_x(\R_+\times \{x\in \Omega: |x|<2R\})}\Bigr)
\\
&\lesssim \sum_{|\alpha|\le 2}\Bigl( \|Z^\alpha f\|_{\dot H^\gamma_B}
+\|Z^\alpha g\|_{\dot H^{\gamma-1}_B}\Bigr)
\notag
\\
&\quad+\sum_{|\alpha|\le 2}
\Bigl(\|\, |x|^{-\frac{n}2+1-\gamma}\Gamma^\alpha F
\|_{L^1_tL^1_rL^2_\omega(\R_+\times \{|x|>2R\})}+\|\Gamma^\alpha F
\|_{L^1_tL_x^{s'_{1-\gamma}}(\R_+\times \{x\in \Omega: |x|<2R\})}\Bigr).\notag
%\\
%&\qquad+\sum_{|\alpha|\le 1}\|\Gamma^\alpha F\|_{L^p_tL_x^{s'_{1-\gamma}}
%(\R_+\times \{x\in \Omega: \, |x|<2R\})}.\notag
\end{align}
Since the $\Gamma$ commute with $\square_g$ when $|x|\ge R$, we have
$$
\square_g \bigl((1-\beta_0)\Gamma^\alpha u\bigr) = 
(1-\beta_0)\Gamma^\alpha F-[\beta_0,\Delta_\g]\Gamma^\alpha u\,.
$$
We can therefore write $(1-\beta_0)\Gamma^\alpha u$ as $v+w$ where 
$\square_g v= (1-\beta_0)\Gamma^\alpha F$ and $v$ has initial data 
$\bigl((1-\beta_0)\Gamma^\alpha u(0,\cd), \partial_t (1-\beta_0)
\Gamma^\alpha u(0,\cd)\bigr)$, while $\square_g w
=-[\beta_0,\Delta_\g]\Gamma^\alpha u$ and $w$ has vanishing initial data.  
If we do this, it follows by \eqref{3.13} that if for $|\alpha|\le 2$ 
we replace the term involving  $(1-\beta_0)\Gamma^\alpha u$ by $v$ in the 
left side of \eqref{a.18}, then the resulting expression is dominated by 
the right side of \eqref{a.18}.  If we use \eqref{14}, we find that if we 
replace $(1-\beta_0)\Gamma^\alpha u$ by $w$ then the resulting expression 
is dominated by
\begin{equation}\label{a.18'}
\sum_{|\alpha|\le 2}
\|\, [\beta_0,\Delta_\g] \Gamma^\alpha u\|_{L^2_tH^{\gamma-1}_B}
\lesssim \sum_{j\le 2}\|\beta_1\partial^j_t u\|_{L^2_tH^{\gamma+2-j}_B}
\end{equation}
assuming that $\beta_1$ equals one on $\supp (\beta_0)$ and is supported 
in  $|x|< 2R$.  As a result, we would be done with the proof of \eqref{a.18} 
if we could show that the right hand side of \eqref{a.18'}
is dominated by the right side of 
\eqref{a.18}. By \eqref{a.14} we control
$\|\beta_1\partial_t^2 u\|_{L^2_tH^{\gamma}_B}$ by the right hand
side of \eqref{a.18}. On the other hand,
$$
\|\beta_1\partial_t u\|^2_{L^2_tH^{\gamma+1}_B}
\lesssim
\|\beta_1\partial^2_t u\|_{L^2_tH^{\gamma}_B}\,
\|\beta_1 u\|_{L^2_tH^{\gamma+2}_B}
$$
so it suffices to dominate $\|\beta_1 u\|_{L^2_tH^{\gamma+2}_B}$.
Since $\Delta_\g u=\partial_t^2u-F$, then 
if $\beta_2$ equals one on $\supp (\beta_1)$ and is supported in the set where 
$|x|<2R$, we may use elliptic regularity and the equation to bound
\begin{align*}
\|\beta_1 u\|_{L^2_tH^{\gamma+2}_B}&\lesssim
\|\beta_2\Delta_\g u\|_{L^2_tH^{\gamma}_B}+
\|\beta_2 u\|_{L^2_tH^{\gamma}_B}\\
&
\lesssim
\|\beta_2\partial_t^2 u\|_{L^2_tH^{\gamma}_B}+
\|\beta_2 u\|_{L^2_tH^{\gamma}_B}+
\|\beta_2 F\|_{L^2_tH^{\gamma}_B}.
\end{align*}
The first two terms are dominated as above using \eqref{a.14}. 
On the other hand, Sobolev embedding and duality yields
\begin{align}\label{a.19}
\|\beta_2 F\|_{L^2_tH^{\gamma}_B}
&\lesssim\sum_{|\alpha|\le 1}\|\partial_x^\alpha F
\|_{L^2_tL^{s'_{1-\gamma}}(\R_+\times\{x\in\Omega:|x|\le 2R\})}\\
&\lesssim\sum_{|\alpha|\le 2}\|\partial_{t,x}^\alpha F
\|_{L^1_tL^{s'_{1-\gamma}}(\R_+\times\{x\in\Omega:|x|\le 2R\})}.
\notag
\end{align}

To finish the proof of \eqref{a.15}, we need to show that the analog of 
\eqref{a.18} is valid when $(1-\beta_0)$ is replaced by $\beta_0$.  Since 
the coefficients of $\Gamma$ are bounded on $\supp(\beta_0)$, if $\beta_1$ 
equals one on $\supp (\beta_0)$ and is supported in $|x|<2R$, then by
Sobolev embedding
\begin{align*}
\sum_{|\alpha|\le 2}\;
\|\beta_0\Gamma^\alpha u
\|_{L^p_tL^{s_\gamma}_x(\R_+\times\Omega)}
&\lesssim\sum_{j\le 2}\;\|\beta_1\partial^j_t u\|_{L^p_tH^{\gamma+2-j}_B}
\\
&\lesssim\sum_{j\le 2}\;\Bigl(\|\beta_1\partial^j_t u\|_{L^2_tH^{\gamma+2-j}_B}
+\|\beta_1\partial^j_t u\|_{L^\infty_tH^{\gamma+2-j}_B}\Bigr).
\end{align*}
The terms in $L^2_tH^{\gamma+2-j}$ are dominated as above.
To control the $L^\infty_tH^{\gamma+2-j}$ terms, and conclude the
proof of \eqref{a.15}, we establish the following estimate:
\begin{multline}\label{a.16}
\sum_{|\alpha|\le 2}
\|\Gamma^\alpha u\|_{L^\infty_t\dot{H}^{\gamma}_B}+
\|\partial_t \Gamma^\alpha u\|_{L^\infty_t\dot{H}^{\gamma-1}_B}
\lesssim \sum_{|\alpha|\le 2}
\Bigl( \|Z^\alpha f\|_{\dot H^\gamma_B(\Omega)}
+\|Z^\alpha g\|_{\dot H^{\gamma-1}_B(\Omega)}\Bigr)
\\
\quad+\sum_{|\alpha|\le 2}
\Bigl(\|\, |x|^{-\frac{n}2+(1-\gamma)}
\Gamma^\alpha F\|_{L^1_tL^1_rL^2_\omega(\R_+\times \{|x|>2R\})}+
\|\Gamma^\alpha F\|_{L^1_tL_x^{s'_{1-\gamma}}
(\R_+\times \{x\in \Omega: |x|<2R\})}\Bigr).
\end{multline}

The inequality where $\Gamma^\alpha u$ is replaced by
$(1-\beta_0)\Gamma^\alpha u$ in \eqref{a.16}
follows by energy estimates on $\R^n$, since the right
hand side dominates $\|(1-\beta_0)\Gamma^\alpha F\|_{L^1_t\dot H^{\gamma-1}}$,
together with \eqref{22'} using the bound \eqref{a.18'} 
to handle the commutator term.
If $\Gamma^\alpha u$ is replaced on the 
left hand side by $\beta_0\Gamma^\alpha u$, 
the result is dominated by 
$\sum_{j\le 3}\|\beta_1\partial_t^j u\|_{L^\infty_tH^{2+\gamma-j}_B}\,.$
For the case $j=0,1$, we write 
$\Box_\g(\beta_1 u)=\beta_1 F-[\Delta_\g,\beta_1] u$, 
and use \eqref{22} with the Duhamel formula to bound
\begin{multline*}
\|\beta_1 u\|_{L^\infty_tH^{\gamma+2}_B}+
\|\beta_1 \partial_t u\|_{L^\infty_tH^{\gamma+1}_B}\\
\lesssim
\|\beta_1 f\|_{H^{\gamma+2}_B}+
\|\beta_1 g\|_{H^{\gamma+1}_B}+
\|\beta_2 u\|_{L^2_tH^{\gamma+2}_B}+
\|\beta_1 F\|_{L^1_tH^{\gamma+1}_B}.
\end{multline*}
The term on the right involving $u$ is controlled previously; on the 
other hand, since $F$ satisfies the $H^{\gamma+1}_B$ boundary
conditions, then
$$
\|\beta_1 F\|_{L^1_tH^{\gamma+1}_B}\lesssim\sum_{|\alpha|\le 2}
\|\partial^\alpha_x F\|_{L^1_tL_x^{s'_{1-\gamma}}}\,.
$$
To handle the terms for $j=2,3$ we use the equation to bound
$$
\sum_{j=2,3}\|\beta_1\partial_t^j u\|_{L^\infty_tH^{2+\gamma-j}_B}\le
\sum_{j=0,1}\Bigl(\|\beta_1\partial_t^j\Delta_\g u\|_{L^\infty_tH^{\gamma-j}_B}+
\|\beta_1\partial_t^j F\|_{L^\infty_tH^{\gamma-j}_B}\Bigr).
$$
The terms involving $\Delta_\g u$ are dominated by
$\|\beta_2\partial_t^j u\|_{L^\infty_tH^{\gamma+2-j}_B}$ with $j=0,1$. 
The terms involving
$F$ are controlled for $j=1$ by \eqref{a.15'}, and for $j=0$
by observing that \eqref{a.19} holds with $L^2_t$
replaced by $L^\infty_t$. This completes the proof of \eqref{a.15} and
\eqref{a.16}.\qed

We shall now use these estimates to prove Theorem \ref{exist}.
  
\noindent{\bf Proof of Theorem \ref{exist}:}
We assume  Cauchy data $(f,g)$ satsifying the smallness condition \eqref{1.7},
and let $u_0$ solve the Cauchy problem \eqref{1} with $F=0$.
We iteratively define $u_k$, for $k\ge 1$, by solving
$$
\begin{cases}
(\partial^2_t-\Delta_\g)u_k(t,x)=F_p(u_{k-1}(t,x))\,, 
\quad (t,x)\in \R_+\times \Omega
\\
u(0,\cd)=f %\in \dot H^\gamma_D(\Omega)
\\
\partial_t u(0,\cd)=g%\in \dot H^{\gamma-1}_D(\Omega)
\\
(Bu)(t,x)=0,\quad \text{on } \, \R_+\times \partial \Omega.
\end{cases}
$$
Our aim is to show that if the constant $\varepsilon>0$ in \eqref{1.7} is 
small enough, then so is
\begin{multline*}
M_k = \sum_{|\alpha|\le 2}\, 
\Bigl( \, 
\bigl\|\Gamma^\alpha u_k\bigr\|_{L^\infty_t\dot H^\gamma_B(\R_+\times\Omega)}+
\bigl\|\partial_t\Gamma^\alpha u_k
\bigr\|_{L^\infty_t\dot H^{\gamma-1}_B(\R_+\times\Omega)}
\\+
\bigl\|\, |x|^{\frac{n}2 - \frac{n+1}p -\gamma} \Gamma^\alpha u_k
\bigr\|_{L^p_tL^p_rL^2_\omega(\R_+\times \{|x|>2R\})}
+\|\Gamma^\alpha u_k
\|_{L^p_tL^{s_\gamma}_x(\R_+\times \{x\in \Omega: \, |x|<2R\})}\, \Bigr)
%+\|\Gamma^\alpha u_k\|_{L^\infty_tL^{s_\gamma}_x(\R_+\times \Omega)}
\end{multline*}
for every $k=0,1,2,\dots$ 

For $k=0$, it follows by \eqref{a.15} and \eqref{a.16} that
$M_0\le C_0\varepsilon$, with $C_0$ a fixed constant.
More generally, \eqref{a.15} and \eqref{a.16} yield that
\begin{align}\label{3.18}
M_k\le C_0\varepsilon +C_0\sum_{|\alpha|\le 2} \, 
\Bigl( \,& \bigl\| \, |x|^{-\frac{n}2+1-\gamma}
\Gamma^\alpha F_p(u_{k-1})
\bigr\|_{L^1_tL^1_rL^2_\omega(\R_+\times \{|x|>2R\})}
\\
&+\|\Gamma^\alpha F_p(u_{k-1})
\|_{L^1_tL^{s_{1-\gamma}'}_x(\R_+\times \{x\in \Omega: \, |x|<2R\})}\Bigr)\,.
\notag
\end{align}
Note that our assumption \eqref{1.4} on the nonlinear term $F_p$ implies that 
for small $v$
$$
\sum_{|\alpha|\le 2}|\Gamma^\alpha F_p(v)|\lesssim |v|^{p-1}
\sum_{|\alpha|\le 2}|\Gamma^\alpha v|+|v|^{p-2}
\sum_{|\alpha|\le 1}|\Gamma^\alpha v|^2\,.
$$
Furthermore, since $u_k$ will be locally of regularity 
$H_B^{\gamma+2}\subset L^\infty$ and
$F_p$ vanishes at $0$, it follows that $F_p(u_k)$ satisfies the
$B$ boundary conditions if $u_k$ does.

Since the collection $\Gamma$ contains vectors spanning the tangent space
to $S^{n-1}$, by Sobolev embedding for $n=3,4$ we have
$$
\|v(r\cd)\|_{L^\infty_\omega}
+\sum_{|\alpha|\le 1}\|\Gamma^\alpha v(r\cd)\|_{L^4_\omega}
\lesssim \sum_{|\alpha|\le 2}
\|\Gamma^\alpha v(r\cd)\|_{L^2_\omega}\,.
$$
Consequently, for fixed $t, r>0$
$$
\sum_{|\alpha|\le 2}\|\Gamma^\alpha F_p(u_{k-1}(t,r\cd) )
\|_{L^2_\omega}\lesssim \sum_{|\alpha|\le 2} \|\Gamma^\alpha u_{k-1}(t,r\cd)
\|^p_{L^2_\omega}\,.
$$
By \eqref{3.10}, the first 
summand in the right side of \eqref{3.18} is dominated by 
$C_1M_{k-1}^p\,.$

We next observe that, since $s_\gamma>2$ and $n\le 4$,
it follows by Sobolev embedding on $\{\Omega\cap|x|<2R\}$ that
$$
\|v\|_{L^\infty(x\in\Omega:|x|<2R)}
+\sum_{|\alpha|\le 1}\|\Gamma^\alpha v\|_{L^4(x\in\Omega:|x|<2R)}
\lesssim \sum_{|\alpha|\le 2}
\|\Gamma^\alpha v\|_{L^{s_\gamma}(x\in\Omega:|x|<2R)}\,.
$$
Since $s_{1-\gamma}'<2$, it holds for each fixed $t$ that
$$
\sum_{|\alpha|\le 2}\|\Gamma^\alpha F_p(u_{k-1}(t,\cd))
\|_{L^{s'_{1-\gamma}}(x\in\Omega:|x|<2R)}
\lesssim \sum_{|\alpha|\le 2} \|\Gamma^\alpha u_{k-1}(t,\cd)
\|^p_{L^{s_\gamma}(x\in\Omega:|x|<2R)}\,.
$$
The second
summand in the right side of \eqref{3.18} is thus also dominated by 
$C_1M_{k-1}^p\,,$ and we conclude that $M_k\le C_0\epsilon+2C_0\,C_1M_{k-1}^p$.
For $\epsilon$ sufficiently small, then
\begin{equation}\label{3.19} M_k\le 2\,C_0\varepsilon, \quad k=1,2,3,\dots
\end{equation}

To finish the proof of Theorem \ref{exist} we need to show that 
$u_k$ converges to a solution of the equation 
\eqref{1.3}.  For this it suffices to show that
\begin{multline*}A_k=
\bigl\| \, |x|^{\frac{n}2-\frac{n+1}p-\gamma}(u_k-u_{k-1})\, 
\bigr\|_{L^p_tL^p_rL^2_\omega
(\R_+\times \{|x|>2R\})} 
\\
+\|u_k-u_{k-1}\|_{L^p_tL^{s_\gamma}_x(\R_+\times \{x\in \Omega: \, |x|<2R\})}
\end{multline*}
tends geometrically to zero as $k\to \infty$.  Since 
$|F_p(v)-F_p(w)|\lesssim |v-w|(\, |v|^{p-1}+|w|^{p-1}\, )$ when $v$ and 
$w$ are small, the proof of \eqref{3.19} can be adapted to show that, 
for small $\varepsilon>0$, there is a uniform constant $C$ so that
$$A_k\le CA_{k-1}(M_{k-1}+M_{k-2})^{p-1},$$
which, by \eqref{3.19}, implies that $A_k\le \tfrac12A_{k-1}$ for small 
$\varepsilon$.  Since $A_1$ is finite, the claim follows, which finishes 
the proof of Theorem \ref{exist}.
\qed


\begin{thebibliography}{MA}
\bibitem{Ben} M. Ben-Artzi: 
{\em Regularity and smoothing for some equations of evolution}, in 
"Nonliner Partial Differential Equations and Applications" 
(H. Brezis and J. L. Lions, eds.), London, 1994, pp. 1--12.

\bibitem{BenK} M. Ben-Artzi and S. Klainerman: 
{\em Decay and regularity for the Schr\"odinger equation}, 
J. Anal. Math. {\bf 58} (1992), 25--37.

\bibitem{BSS} M. Blair, H. Smith and C. D. Sogge: 
{\em Strichartz estimates for the wave equation on manifolds with boundary}, 
arXiv:0805.4733.

\bibitem{burq} N. Burq: {\em Global Strichartz estimates for 
nontrapping geometries: About an article by H. Smith and C. Sogge}, 
Comm. Partial Differential Equations  {\bf 28} (2003), 1675--1683.

\bibitem{BLP} N. Burq, G. Lebeau and F. Planchon: 
{\em Global existence for energy critical waves in 3-D domains}, 
J. Amer. Math. Soc. {\bf 21} (2008), 831--845.

\bibitem{BP} N. Burq and F. Planchon: 
{\em Global existence for energy critical waves in 3-d domains : 
Neumann boundary conditions}, arXiv:0711.0275.

\bibitem{CK} M. Christ and A. Kiselev: 
{\em Maximal functions associated to filtrations}, 
J. Funct. Anal. {\bf 179} (2001), 409--425.

\bibitem{MSYY}Y. Du, J. Metcalfe, C. D. Sogge and Y. Zhou: 
{\em Concerning the Strauss conjecture and almost global existence 
for nonlinear Dirichlet-wave equations in $4$-dimensions},
Comm. Partial Differential Equations {\bf 33} (2008), 1487--1506.

\bibitem{YY} Y. Du and Y. Zhou: 
{\em The life span for nonlinear wave equation outside of star-shaped
obstacle in three space dimensions}, Comm. Partial Differential
Equations {\bf 33} (2008), 1455--1486.

\bibitem{FW} D. Fang and C. Wang: 
{\em Weighted Strichartz Estimates with Angular Regularity and their 
Applications}, arXiv:0802.0058.

\bibitem{GLS} V. Georgiev, H. Lindblad, and C. D. Sogge: 
{\em Weighted Strichartz estimates and global existence for 
semilinear wave equations},  
Amer. J. Math.  {\bf 119}  (1997),  1291--1319.

\bibitem{Glassey} R. T. Glassey: 
{\em Existence in the large for $\square u = F(u)$ in two dimensions}, 
Math. Z. {\bf 178} (1981), 233--261.

\bibitem{H} K. Hidano: 
{\em Morawetz-Strichartz estimates for spherically symmetric solutions 
to wave equations and applications to semi-linear Cauchy problems}, 
Differential Integral Equations {\bf 20} (2007), 735--754.

\bibitem{H2}  K. Hidano: 
{\em Small solutions to semi-linear wave equations with radial data 
of critical regularity}, Rev. Mat. Iberoamericana, to appear.

\bibitem{Hos}{T. Hoshiro}: 
{\em On weighted $L^2$ estimates of solutions to wave equations}, 
J. Anal. Math. {\bf 72} (1997), 127--140.

\bibitem{John} F. John: 
{\em Blow-up of solutions of nonlinear wave equations in three space 
dimensions},  
Manuscripta Math. {\bf 28} (1979), 235--265

\bibitem{KSS}  M. Keel, H. F. Smith, and C. D. Sogge: 
{\em Almost global existence for some semilinear wave equations}, 
J. Anal. Math. {\bf 87} (2002), 265-279.

\bibitem{KT} M. Keel and T. Tao, 
{\em Endpoint Strichartz estimates}, 
Amer. J. Math. {\bf 120} (1998), 955--980.

\bibitem{levine} H. A. Levine, {\em Instability and nonexistence of global solutions to nonlinear wave equations of the form} $Pu\sb{tt}=-Au+{\cal F}(u)$, Trans. Amer. Math. Soc. {\bf 192} (1974), 1--21.

\bibitem{LX}{T. T. Li and Y. Zhou}: 
{\em A note on the life-span of classical solutions to nonlinear wave 
equations in four space dimensions}, 
Indiana Univ. Math. J. {\bf 44} (1995), 1207--1248.

\bibitem{LS}{H. Lindblad and C. D. Sogge}: 
{\em On existence and scattering with minimal regularity for semilinear 
wave equations}, 
J. Funct. Anal. {\bf 130} (1995), 357--426.

\bibitem{LS2}{H. Lindblad and C. D. Sogge}: 
{\em Long-time existence for small amplitude semilinear wave equations}, 
Amer. J. Math. {\bf 118} (1996), 1047--1135.

\bibitem{Mel} R. B. Melrose: 
{\em Singularities and energy decay in acoustical scattering}, 
Duke Math. J. {\bf 46} (1979), 43--59.

\bibitem{MelSj} R. B. Melrose and J. Sj\"ostrand: 
{\em Singularities of boundary value problems. I}, 
Comm. Pure Appl. Math. {\bf 31} (1978), 593--617.

\bibitem{M}{J. Metcalfe}: 
{\em Global Strichartz estimates for solutions to the wave equation 
exterior to a convex obstacle},   
Trans. Amer. Math. Soc.  {\bf 356}, (2004), 4839--4855.

\bibitem{mor} C. S. Morawetz: 
{\em Decay for solutions of the exterior problem for the wave equation}, 
Comm. Pure and Appl. Math.  {\bf 28} (1975), 229--264.

\bibitem{mrs} C. S. Morawetz, J. Ralston and W. Strauss: 
{\em Decay of solutions of the wave equation outside nontrapping obstacles},
Comm. Pure Appl. Math. {\bf 30} (1977), 87--133.

\bibitem{Lax} P. D. Lax and R. S. Philips: 
Scattering Theory (Revised Edition), Academic Press Inc., 1989.

\bibitem{Ral} J. Ralston: 
{\em Note on the decay of acoustic waves}, 
Duke Math. J. {\bf 46} (1979), 799-804.

\bibitem{ST} Y. Shibata and Y. Tsutsumi: 
{\em Global existence theorem for nonlinear wave equation in exterior domain},
Lecture Notes in Num. Appl. Anal., Vol. 6 (1983), 155-196.

\bibitem{Sideris} T. C. Sideris: 
{\em Nonexistence of global solutions to semilinear wave equations 
in high dimensions}, 
J. Differential Equations {\bf 52} (1984), 378--406.

\bibitem{SS0} H. F. Smith and C. D. Sogge: 
{\em On the critical semilinear wave equation outside convex obstacles}, 
J. Amer. Math. Soc. {\bf 8} (1995), 879--916.

\bibitem{SS} H. F. Smith and C. D. Sogge: 
{\em Global Strichartz estimates for nontrapping perturbations of 
the Laplacian}, 
Comm. Partial Differential Equations {\bf 25}, (2000), 2171--2183.

\bibitem{SS2} H. F. Smith and C. D. Sogge: 
{\em On the $L\sp p$ norm of spectral clusters for compact manifolds 
with boundary}, Acta Math. {\bf 198}, (2007), 107--153.

\bibitem{So} C. D. Sogge: 
Lectures on nonlinear wave equations, International Press, Boston, MA 1995.

\bibitem{So2} C. D. Sogge: 
Lectures on nonlinear wave equations, 2nd edition, 
International Press, Boston, MA, 2008.

\bibitem{Ta} D. Tataru: 
{\em Strichartz estimates in the hyperbolic space and global existence 
for the semilinear wave equation},  
Trans. Amer. Math. Soc.  {\bf 353}  (2001),   795--807.

\bibitem{taylor} M. Taylor:
{\em Grazing rays and reflection of singularities of solutions
to wave equations},
Comm. Pure Appl. Math. {\bf 29} (1976), 1--38.

\bibitem{V} B. R. Vainberg: 
{\em The short-wave asymptotic behavior of the solutions of stationary 
problems, and the asymptotic behavior as $t\rightarrow \infty $ of the 
solutions of nonstationary problems}, 
Russian Math. Surveys {\bf 30} (1975), 1--58.

\bibitem{Zhou} Y. Zhou: 
{\em Cauchy problem for semilinear wave equations with small data in 
four space dimensions},
J. Partial Differential Equations {\bf 8} (1995), 135--144.


\end{thebibliography}
\end{document}